\documentclass[10pt, a4paper, reqno]{amsart}

\usepackage[english]{babel}
\usepackage{amssymb}
\usepackage{amsmath}
\usepackage{enumerate}
\usepackage{graphicx}
\usepackage[all]{xy}
\usepackage{longtable}
\usepackage{multirow}
\usepackage{cancel}
\usepackage{hyperref}

\newtheorem{Thm}{Theorem}[section]
\newtheorem{Lemma}[Thm]{Lemma}
\newtheorem{Prop}[Thm]{Proposition}

\newtheorem*{MQuest}{Main Question}
\newtheorem*{RecallThm}{Main Theorem}

\theoremstyle{definition}
\newtheorem{Def}[Thm]{Definition}
\newtheorem{Ex}[Thm]{Example}

\theoremstyle{remark}
\newtheorem{Rmk}[Thm]{Remark}

\numberwithin{equation}{section}

\newcommand{\C}{\mathbb{C}}
\newcommand{\R}{\mathbb{R}}
\newcommand{\Z}{\mathbb{Z}}
\newcommand{\N}{\mathbb{N}}

\newcommand{\F}{\mathbb{F}}
\newcommand{\PP}{\mathbb{P}}
\newcommand{\cO}{\mathcal{O}}
\newcommand{\cL}{\mathcal{L}}
\newcommand{\cT}{\mathcal{T}}
\newcommand{\cN}{\mathcal{N}}
\newcommand{\cF}{\mathcal{F}}

\newcommand{\chitop}{\chi_{\text{top}}}

\newcommand{\Pic}{\operatorname{Pic}}
\newcommand{\MW}{\operatorname{MW}}

\newcommand{\id}{\operatorname{id}}

\title{Families of Calabi--Yau elliptic fibrations in \texorpdfstring{$\PP(\cL^\MakeLowercase{a} \oplus \cL^\MakeLowercase{b} \oplus \cO_B)$}{P(La + Lb + O)}}

\author{Andrea Cattaneo}
\address{Andrea Cattaneo, Dipartimento di Scienza e Alta Tecnologia, Universit\`a degli Studi dell'Insubria, Via Valleggio 11, 22100, Como, Italy}

\address{Andrea Cattaneo, Institut Camille Jordan UMR 5208, Universit\'e Claude Bernard Lyon 1, 69622 Villeurbanne Cedex, France}
\email{cattaneo@math.univ-lyon1.fr}

\begin{document}

\begin{abstract}
Let $B$ be a smooth projective surface, and $\cL$ an ample line bundle on $B$. The aim of this parer is to study the families of elliptic Calabi--Yau threefolds sitting in the bundle $\PP(\cL^a \oplus \cL^b \oplus \cO_B)$ as anticanonical divisors. We will show that the number of such families is finite.
\end{abstract}

\thanks{This paper collects part of the author's Ph.D.~thesis: he wants to gratefully acknowledge his advisor, Prof.~Bert van Geemen, for the many suggestions and his support during the preparation and the writing of the paper. He wants also to thanks the referee for his/her useful advices, especially for having focussed the attention on the case of Hirzebruch surfaces.}

\subjclass[2010]{Primary 14J30, 14J32.}

\keywords{Elliptic threefolds, Calabi--Yau varieties.}

\maketitle

\tableofcontents

\section*{Introduction}
While the theory of elliptic surfaces is a well settled and consolidated subject, in the case of elliptic threefolds there are still many interesting and open questions. Not only the theoretical aspects of the theory are important, but also the research of families of examples plays a central role: one of the main motivations is the close connection with the theory of strings (and in particular $F$-theory, see e.g.~\cite{EvidenceForFTheory}), which is a physical subject whose main object of study is in fact elliptic fibrations on Calabi--Yau manifolds. To give two examples, in \cite{EsoleYau} the $E_6$ and $E_7$ family of elliptic Calabi--Yau threefolds are defined, and in \cite{EsoleYauD5} the authors define the $D_5$ family.

In this paper we will focus on a way of constructing elliptic fibrations on Calabi--Yau threefold, working on the field $\mathbb{C}$ of complex numbers.

A simple way to produce Calabi--Yau varieties is to consider smooth anticanonical subvarieties of some reasonable ambient space: in fact by adjunction these varieties will automatically be Calabi--Yau. Giving different shades to the word `reasonable', one has different classes of ambient spaces to try describing its anticanonical subvarieties. In particular the class of toric Fano Gorenstein fourfolds has been deeply studied for the following reasons:
\begin{enumerate}
\item Since any anticanonical divisor of a Fano variety is ample, we are sure to find effective divisors in the anticanonical system;
\item Gorenstein varieties may be singular, but in this case they have nice resolutions of the singularities and one can then study the anticanonical subvarieties of the resolution;
\item Toric varieties are simple since most of the problems one may have to solve can be translated into a combinatorial problem, which is simpler to deal with.
\end{enumerate}
To each toric Fano Gorenstein fourfold is associated a reflexive $4$-dimensional polyhedron and viceversa, so the first attempt to describe the Calabi--Yau subvarieties in these ambient space is to classify all the reflexive $4$-dimensional polyhedra. Such a classification is known, and there are $473,800,776$ $4$-dimensional reflexive polyhedra (see e.g.~\cite{KreuzerSkarke}, \cite{KreuzerSkarke2}). Among these, in \cite{Braun} the $102,581$ elliptic fibrations over $\PP^2$ are identified.

The elliptic fibrations we will describe in this paper are anticanonical hypersurfaces in a projective bundle $Z$ over a surface $B$ of the form $Z = \PP(\cL^a \otimes \cL^b \otimes \cO_B)$ for $\cL$ an ample line bundle on $B$. Observe that even in the case where the base $B$ is toric, e.g.~$B = \PP^2$, the ambient bundle is typically not Fano.

The aim of this paper is to show that once $B$ and $\cL$ are fixed, then the bundle $Z = \PP(\cL^a \oplus \cL^b \oplus \cO_B)$ can house Calabi--Yau elliptic fibrations only for a finite number of choices of $(a, b)$:
\begin{RecallThm}[Thm. \ref{thm: finiteness result}]
Let $B$ be a smooth projective surface, and $\cL$ an ample line bundle on $B$. Consider the projective bundle $\PP(\cL^a \oplus \cL^b \oplus \cO_B)$, with $a \geq b \geq 0$. Then only for a finite number of pairs $(a, b)$ the generic anticanonical hypersurface in $\PP(\cL^a \oplus \cL^b \oplus \cO_B)$ is a Calabi--Yau elliptic fibration over $B$.
\end{RecallThm}
As we will see in sections~\ref{section: statement of the problem} and \ref{section: proof theorem}, we may fail to find a Calabi--Yau elliptic fibration for the following reasons: the fibration has no sections or its total space is singular.

The outline of the paper is as follows. In section~\ref{section: definitions} we will recall the definitions of elliptic fibration and of Calabi--Yau variety. In section~\ref{section: the result} we will state the finiteness result (theorem~\ref{thm: finiteness result}), and prove it (sections from \ref{section: step1} to \ref{section: reducible conic case}). Finally, in section~\ref{section: examples} we will give some concrete examples: we will find explicit bounds on the number of different families when the base $B$ is a del Pezzo surface (and in particular for $B = \PP^2$), and when the basis $B$ is a Hirzebruch surface $\F_e$.

\section{Elliptic fibrations and Calabi--Yau manifolds}\label{section: definitions}
In this section, we want to recall the definition and main properties of elliptic fibrations (section~\ref{section: elliptic fibrations}) and Calabi--Yau manifolds (section~\ref{section: calabi-yau varieties}). Throughout this paper, all the varieties are defined over $\C$.

\subsection{Elliptic fibrations}\label{section: elliptic fibrations}
Elliptic fibrations are the geometric realization of elliptic curves over the function field of a variety. Their study has been encouraged by physics, and in particular string theory: to each elliptic fibration correspond a physical scenario, and the fibration itself determines the number of elementary particles, their charges and masses (see e.g.~\cite{EvidenceForFTheory}).

\begin{Def}
We say that $\pi: X \longrightarrow B$ is an \emph{elliptic fibration} over $B$ if
\begin{enumerate}
\item $X$ and $B$ are projective varieties of dimension $n$ and $n - 1$ respectively, with $X$ smooth;
\item $\pi$ is a surjective morphism with connected fibres;
\item the generic fibre of $\pi$ is a smooth connected curve of genus $1$;
\item a section $\sigma: B \longrightarrow X$ of $\pi$ is given.
\end{enumerate}
When $\pi: X \longrightarrow B$ satisfies only the first three requirements above, we say that it is a \emph{genus one fibration}.
\end{Def}

We will denote the fibre over the point $P \in B$ with $X_P$.

\begin{Rmk}
Let $\pi: X \longrightarrow B$ be an elliptic fibration, with section $\sigma$. Then each smooth fibre $X_P$ is an elliptic curve, where we choose as origin the point $\sigma(P)$.
\end{Rmk}

A \emph{morphism} between two elliptic fibrations $\pi: X \longrightarrow B$ and $\pi': X' \longrightarrow B$ is a morphism of varieties over $B$, i.e.\ a morphism $f: X \longrightarrow X'$ such that
\[\xymatrix{X \ar[rr]^f \ar[dr]_\pi& & X' \ar[dl]^{\pi'}\\
 & B & }\]
commutes.

Not every fibre of $\pi$ needs to be smooth: the \emph{discriminant locus} of the fibration is the subset of $B$ over which the fibres are singular
\[\Delta = \{ P \in B \,|\, X_P \text{ is singular} \} \subseteq B.\]

A \emph{rational section} of $\pi$ is a rational map $s: B \dashrightarrow X$ such that $\pi \circ s = \id$ over the domain of $s$. The \emph{Mordell--Weil group} of the fibration is
\[\MW(X) = \{ s: B \dashrightarrow X \,|\, s \text{ is a rational section} \},\]
where the group law is given by addition fibrewise. Observe that even though the elements of the Mordell--Weil group are rational sections, we require its zero element to be a section.

\subsubsection{The Weierstrass model of an elliptic fibration}\label{section: Weierstrass model}

The main reason for requiring that an elliptic fibration admits a section is that we can use the presence of this section to define the \emph{Weierstrass model} of the fibration (see \cite[Thm.~2.1]{OnWeiMod}).

Let $\pi: X \longrightarrow B$ be an elliptic fibration. By a little abuse of notation, we still call the image of the distinguished section, $S = \sigma(B)$, the distinguished section of $X$. Denote by $i$ the inclusion $i: S \hookrightarrow X$, then we define the \emph{fundamental line bundle} of the fibration as the line bundle on $B$
\[\cF = \left( \pi_* i_* \cN_{S|X} \right)^{-1},\]
and the \emph{Weierstrass model} of $X$ is then the image of the birational morphism
\[f: X \longrightarrow \PP(\pi_* \cO_X(3 S)) = \PP(\cF^{\otimes 2} \oplus \cF^{\otimes 3} \oplus \cO_B)\]
defined by the canonical morphism $\pi^* \pi_* \cO_X(3S) \longrightarrow \cO_X(3S)$. For the surjectivity of this map we refer to \cite[Proof of Thm.~2.1]{OnWeiMod}.

\begin{Rmk}
Let $p: W \longrightarrow B$ be the Weierstrass model of $\pi: X \longrightarrow B$. Then $W$ is defined in $\PP(\cF^{\otimes 2} \oplus \cF^{\otimes 3} \oplus \cO_B)$ by a Weierstrass equation
\begin{equation}\label{eq: Weierstrass eq}
W: y^2 z = x^3 + \alpha_{102} x z^2 + \alpha_{003} z^3,
\end{equation}
where $\alpha_{102} \in H^0(B, \cF^{\otimes 4})$, $\alpha_{003} \in H^0(B, \cF^{\otimes 6})$.
\end{Rmk}

\begin{Rmk}
The discriminant locus $\Delta$ of a Weierstrass fibration $p: W \longrightarrow B$ is not only a subset of $B$, but also a subvariety (actually, a subscheme) of $B$. It is defined in terms of the coefficients of the Weierstrass equation \eqref{eq: Weierstrass eq} by
\[\Delta: 4 \alpha_{102}^3 + 27 \alpha_{003}^2 = 0.\]
\end{Rmk}

\subsection{Calabi--Yau manifolds}\label{section: calabi-yau varieties}
Calabi--Yau manifolds are the higher dimensional analogues of elliptic curves and $K3$ surfaces. The mathematical models of $F$-theory are all examples of Calabi--Yau manifolds: this property is needed on the total space of an elliptic fibration in order to get a physically consistent model (see e.g.~\cite{CompactificationsOfFTheoryI, CompactificationsOfFTheoryII}).

\begin{Def}
A \emph{Calabi--Yau manifold} is a smooth compact K\"ahler variety $X$ with
\begin{enumerate}
\item trivial canonical bundle $\omega_X \simeq \cO_X$,
\item $h^{0, q}(X) = 0$ for $q = 1, \ldots \dim X - 1$, where $h^{p, q}(X) = \dim H^q(X, \Omega^p_X)$.
\end{enumerate}
\end{Def}

\begin{Ex}
If $X$ is a Calabi--Yau variety of dimension $1$, then $X$ is a smooth Riemann surface of genus $1$.

In the case of dimension $2$, the Calabi--Yau surfaces are the $K3$.

In dimension $3$, the Fermat quintic in $\PP^4$, and in fact any smooth quintic, is a classical example of Calabi--Yau variety (see for instance \cite{GHJ} and \cite{CoxKatz}). Other Calabi--Yau threefolds which are complete intersections in projective spaces are the complete intersection of two hypersurfaces of degree $3$ in $\PP^5$, of a hyperquadric and a hypersurface of degree $4$ in $\PP^5$, of two hyperquadric and a hypercubic in $\PP^6$ or the complete intersection of four hyperquadrics in $\PP^7$.

For other examples of Calabi--Yau manifolds, see e.g.~\cite{Bestiary}.
\end{Ex}

\section{A finiteness result}\label{section: the result}

\subsection{Notations and general setting}
In this section we will fix the notation we will use through the rest of the paper.

Let $B$ be a smooth projective surface, and $\pi: X \longrightarrow B$ an elliptic threefold over $B$. As we observed in section~\ref{section: Weierstrass model}, the Weierstrass model of $\pi$ sits in a projective bundle of the form $\PP(\cF^{\otimes 2} \oplus \cF^{\otimes 3} \oplus \cO_B)$ for a suitable line bundle $\cF$ on $B$. We then want to investigate all the elliptic fibrations that can be embedded in similar ambient spaces.

This is the general framework we will be working in. Let $B$ be a smooth projective surface, and $\cL$ an ample line bundle on $B$. Let $p: Z \longrightarrow B$ be the projective bundle \emph{of lines} associate to the rank two vector bundle $\cL^a \oplus \cL^b \oplus \cO_B$, i.e.\ $Z = \PP(\cL^a \oplus \cL^b \oplus \cO_B)$.

Let $X \in |-K_Z|$ be an anticanonical subvariety, and $\pi: X \longrightarrow B$ the restriction to $X$ of the structure map $p$ of $Z$.

\subsection{Statement of the problem}\label{section: statement of the problem}
The aim of the paper is to give an answer to the following question:

\begin{MQuest}
For how many (and for which) pairs $(a, b)$ is it true that for the generic anticanonical subvariety $X$ of $\PP(\cL^a \oplus \cL^b \oplus \cO_B)$, the map $\pi$ defines a Calabi--Yau elliptic fibration over $B$?
\end{MQuest}

At first sight the answer seems to be ``almost for all pairs'', for the following reasons:
\begin{enumerate}
\item anticanonical subvarieties are Calabi--Yau by adjunction;
\item since the generic fibre of $\pi$ is a plane cubic curve (cf.~\eqref{equation: cubic equation}), we have always a genus $1$ fibration.
\end{enumerate}

Nevertheless we are wrong. In fact the map $\pi$ can have no sections, or the total space $X$ of the fibration can be singular. This last case can happen for two reasons:
\begin{enumerate}
\item the generic $X \in |-K_Z|$ is reducible (see section~\ref{section: reducible conic case});
\item there is a section of $\pi$ passing through a singular point of a fibre.
\end{enumerate}
In the second case, if the singularities of $X$ admit a small resolution we can obtain a Calabi--Yau elliptic fibration, but then the resolved fibration would live in another ambient space, so we exclude them from this paper.

\begin{Thm}\label{thm: finiteness result}
Let $B$ be a smooth projective surface, and $\cL$ an ample line bundle on $B$. Consider the projective bundle $\PP(\cL^a \oplus \cL^b \oplus \cO_B)$, with $a \geq b \geq 0$. Then only for a finite number of pairs $(a, b)$ the generic anticanonical hypersurface in $\PP(\cL^a \oplus \cL^b \oplus \cO_B)$ is a Calabi--Yau elliptic fibration over $B$.
\end{Thm}

\begin{Rmk}
Our theorem can be considered as a reflex of the more general statement \cite[Thm.~0.1]{Gross} that there are only a finite number of deformation families of Calabi--Yau elliptic threefolds over rational surfaces with the property that any Calabi--Yau elliptic threefold over a rational surface is birational to one elliptic fibration in these families (see also \cite[Thm.~1.1]{DiCerboSvaldi} for an analogue statement for Calabi--Yau elliptic fourfolds and fivefolds).
\end{Rmk}

\begin{Rmk}
The theorem states only the finiteness, but its proof gives also a sort of algorithm to detect a finite superset of the set of pairs satisfying the main question.
\end{Rmk}

\begin{Rmk}
Consider the projective bundle $\PP(\cL^\alpha \oplus \cL^\beta \oplus \cO_B)$, with $(\alpha, \beta) \in \Z \times \Z$. After tensoring $\cL^\alpha \oplus \cL^\beta \oplus \cO_B$ with $\cL^{-m}$, where $m = \min\{ \alpha, \beta, 0 \}$, and a permutation of the addends, we get a new vector bundle, of the form $\cL^a \oplus \cL^b \oplus \cO_B$ with $a \geq b \geq 0$, and such that
\[\PP(\cL^\alpha \oplus \cL^\beta \oplus \cO_B) \simeq \PP(\cL^a \oplus \cL^b \oplus \cO_B).\]
So the bound on the possible $(a, b)$'s in the hypothesis of theorem~\ref{thm: finiteness result} is not restrictive.
\end{Rmk}

Before proving theorem~\ref{thm: finiteness result}, in section~\ref{section: the ambient bundle} we will take a short digression on the projective bundle $Z$ and its anticanonical subvarieties.

\subsection{Calabi--Yau's in \texorpdfstring{$\PP(\cL^a \oplus \cL^b \oplus \cO_B)$}{P(La + Lb + O)}}\label{section: the ambient bundle}

We are interested in studying the anticanonical subvarieties of $Z = \PP(\cL^a \oplus \cL^b \oplus \cO_B)$. In this section we want first to compute the Chern classes of $Z$, and then find how an equation for an anticanonical subvariety looks like.

\subsubsection{The ambient bundle}\label{section: chern classes of Z}
The bundle projection $p: Z \longrightarrow B$ gives the relative tangent bundle exact sequence
\[0 \longrightarrow \cT_{Z|B} \longrightarrow \cT_Z \longrightarrow p^* \cT_B \longrightarrow 0\]
from which we see that
\[c(Z) = c(\cT_{Z|B}) p^* c(B).\]
To compute the total Chern class of the relative tangent bundle, we exploit the fact that it fits into an Euler-type exact sequence (see \cite[p. 435, B.5.8]{FulInt}):
\[0 \longrightarrow \cO_Z \longrightarrow p^*E \otimes \cO_Z(1) \longrightarrow \cT_{Z|B} \longrightarrow 0,\]
where $E = \cL^a \oplus \cL^b \oplus \cO_B$.

An explicit computation leads to the following results
\begin{equation}\label{formula: chern classes of Z}
\begin{array}{rl}
c_1(Z) = & p^* c_1(B) + (a + b) p^* L + 3 \xi,\\
c_2(Z) = & ab p^* L^2 + (a + b) p^* L c_1(B) + 2(a + b) p^*L \xi +\\
 & + 3 p^* c_1(B) \xi + p^* c_2(B) + 3 \xi^2,\\
c_3(Z) = & 2(a + b) p^* c_1(B) L \xi + 3 p^*c_1(B) \xi^2 + 3 p^* c_2(B) \xi,\\
c_4(Z) = & 3 p^* c_2(B) \xi^2,\\
\end{array}
\end{equation}
where $L = c_1(\cL)$ and $\xi = c_1(\cO_Z(1))$.

\subsubsection{Equations for anticanonical subvarieties}
Consider the projective bundle $Z = \PP(\cL^a \oplus \cL^b \oplus \cO_B)$, and let $x$, $y$ and $z$ denote sections on $Z$ whose vanishing gives the subvariety of $Z$ corresponding to the embeddings
\[\cL^b \oplus \cO_B \hookrightarrow E, \qquad \cL^a \oplus \cO_B \hookrightarrow E, \qquad \cL^a \oplus \cL^b \hookrightarrow E\]
respectively. Then
\[\begin{array}{l}
x \in H^0(Z, p^*\cL^a \otimes \cO_Z(1))\\
y \in H^0(Z, p^*\cL^b \otimes \cO_Z(1))\\
z \in H^0(Z, \cO_Z(1))
\end{array}\]
and we can use $(x: y: z)$ as global homogeneous coordinates in $Z$ over $B$.

Since $c_1(Z) = p^* c_1(B) + (a + b) p^* L + 3 \xi$ by \eqref{formula: chern classes of Z}, an equation $F$ defining an anticanonical hypersurface must be cubic in $(x: y: z)$, of the form
\begin{equation}\label{equation: cubic equation}
F = \sum_{i + j + k = 3} \alpha_{ijk} x^i y^j z^k,
\end{equation}
and the coefficient $\alpha_{ijk}$ of the monomial $x^i y^j z^k$ must be a section of a suitable line bundle, according to table~\ref{table: weights}.

\begin{center}
\begin{longtable}{|c|c|}
\caption{Cubic monomials and the weight of their coefficients.}
\label{table: weights}\\
\hline
Monomial & Weight of the coefficient\\
\hline
$x^3$ & $c_1(B) - 2aL + bL$\\
\hline
$x^2 y$ & $c_1(B) - aL$\\
\hline
$x y^2$ & $c_1(B) - bL$\\
\hline
$y^3$ & $c_1(B) + aL - 2bL$\\
\hline
$x^2 z$ & $c_1(B) - aL + bL$\\
\hline
$xyz$ & $c_1(B)$\\
\hline
$y^2 z$ & $c_1(B) + aL - bL$\\
\hline
$x z^2$ & $c_1(B) + bL$\\
\hline
$y z^2$ & $c_1(B) + aL$\\
\hline
$z^3$ & $c_1(B) + aL + bL$\\
\hline
\end{longtable}
\end{center}

\subsubsection{Chern classes of anticanonical subvarieties}\label{ChernChapter}
We want to compute the Chern classes of a smooth $X \in |-K_Z|$. We have
\[\xymatrix{X \ar@{^(->}[rr]^i \ar[dr]_\pi & & Z \ar[dl]^p\\
 & B & }\]
and the normal bundle sequence of $X$ in $Z$
\[0 \longrightarrow \cT_X \longrightarrow i^* \cT_Z \longrightarrow \cN_{X | Z} \longrightarrow 0,\]
which gives the following relation between the total Chern classes
\[i^* c(Z) = c(X) c(\cN_{X | Z}) = c(X) i^*(1 - K_Z).\]

Since we know $c(Z)$ from section~\ref{section: chern classes of Z}, and $1 - K_Z$ is a unit in the Chow ring of $Z$, we deduce the following formulae for the Chern classes of $X$:
\begin{equation}\label{formula: chern classes of X}
\begin{array}{rl}
c_1(X) = & 0,\\
c_2(X) = & 3 \xi_{|_X}^2 + \pi^*(2(a + b)L + 3 c_1(B)) \xi_{|_X} +\\
 & + \pi^*((a + b)L c_1(B) + ab L^2 + c_2(B)),\\
c_3(X) = & - 9 \pi^* c_1(B) \xi_{|_X}^2 - \pi^* (2 (a^2 - ab + b^2) L^2 +\\
 & + 6 (a + b) L c_1(B) + 3 c_1(B)^2) \xi_{|_X}.
\end{array}
\end{equation}

\begin{Rmk}
In particular, we have a formula for the Euler--Poincar\'e characteristic of our varieties:
\[\chitop(X) = \deg c_3(X) = -6(a^2 - ab + b^2)L^2 - 18 c_1(B)^2.\]
\end{Rmk}

\subsection{Hypersurfaces in Calabi--Yau threefolds}

In this section we want to recall a result which will be crucial in the proof of our main theorem~\ref{thm: finiteness result} (see section~\ref{section: step2}). Using proposition~\ref{prop: friedman} we will in fact reduce our general problem to a simpler one, concerning only the base surface of the elliptic fibration.

Assume that $X$ is any threefold with $c_1(X) = 0$, and that $i: S \hookrightarrow X$ is the inclusion of a smooth surface. The techniques used in section~\ref{section: chern classes of Z} and section~\ref{ChernChapter} can be used to have more information on how $S$ is embedded in $X$.

From the normal bundle sequence
\[0 \longrightarrow T_S \longrightarrow i^* T_X \longrightarrow \cN_{S|X} \longrightarrow 0\]
we get that
\begin{equation}\label{formula: chern relation normal bundle}
i^* c(X) = c(S) c(\cN_{S|X}).
\end{equation}
To compute $i_* c(\cN_{S|X})$, we can argue in two ways:
\begin{itemize}
\item By the self-intersection formula, $c(\cN_{S|X}) = i^*(1 + [S])$ where $[S]$ is the class of $S$ in the Chow ring of $X$. So
\begin{equation}\label{formula: nbund 1}
i_* c(\cN_{S|X}) = i_* i^* (1 + [S]) = (1 + [S]) [S] = [S] + [S]^2.
\end{equation}
\item Using \eqref{formula: chern relation normal bundle}, we have that $c(\cN_{S|X}) = i^* c(X) \cdot c(S)^{-1}$, and so
\begin{equation}\label{formula: nbund 2}
i_* c(\cN_{S|X}) = c(X) \cdot i_*(c(S)^{-1}) = [S] - i_* c_1(S) + c_2(X)[S] - i_*(c_2(S) - c_1(S)^2).
\end{equation}
\end{itemize}
Comparing \eqref{formula: nbund 1} and \eqref{formula: nbund 2} we get that
\begin{equation}\label{formula: intersection results}
[S]^2 = -i_* c_1(S), \qquad c_2(X)[S] = i_*(c_2(S) - c_1(S)^2)
\end{equation}

Taking the degree of the second relation in \eqref{formula: intersection results} gives us the following result.

\begin{Prop}[{\cite[Lemma~4.4]{Friedman}}]\label{prop: friedman}
Let $X$ be a threefold with $c_1(X) = 0$, and $S$ a smooth hypersurface, with associated class $[S]$. Then
\[c_2(X)[S] = \chitop(S) - K_S^2.\]
\end{Prop}

The first relation in \eqref{formula: intersection results} gives an interpretation to $[S]^2$. To understand also what $[S]^3$ is we use the adjunction formula for $S$ in $X$:
\[c_1(S) = i^*(c_1(X) - [S]) = -i^*[S].\]
From this relation we have that
\[\deg c_1(S)^2 = \deg i^*[S]^2 = \deg i_* i^*[S]^2 = \deg [S]^3,\]
i.e.~$K_S^2 = [S]^3$.

\subsubsection{The fundamental line bundle of a Calabi--Yau elliptic fibration}\label{sect: fund line bund}

Assume that $\pi: X \longrightarrow B$ is an elliptic fibration with section $S$, where $X$ a Calabi--Yau threefold. We can use the first relation in \eqref{formula: intersection results} to compute the fundamental line bundle of $\pi$. In fact, since $\pi_*([S]) = B$ we have that
\[\pi_* i_* c(\cN_{S|X}) = \pi_*([S] + [S]^2) = 1 - p_* i_* c_1(S) = 1 - c_1(B).\]
So, if $\cF$ is the fundamental line bundle of $\pi$, then $c_1(\cF) = c_1(B)$, and so we can embed the Weierstrass model of $\pi$ in
\[\PP(\omega_B^{-2} \oplus \omega_B^{-3} \oplus \cO_B),\]
where $\omega_B$ is the anticanonical line bundle of $B$.

\subsection{Proof of the Main Theorem}\label{section: proof theorem}
We will split the proof of theorem~\ref{thm: finiteness result} in several steps to make it clearer. In the first step (section~\ref{section: step1}) we will show that possibly with the exception of a finite number of pairs $(a, b)$, the genus one fibrations $X$ in $Z = \PP(\cL^a \oplus \cL^b \oplus \cO_B)$ admit a section. In the second step (section~\ref{section: step2}) we will concentrate on such pairs, and use the presence of the section to reduce the problem to a new problem concerning only the intersection form on the base. In the third step (sections~\ref{section: step 3irr} and \ref{section: reducible conic case}) we will show that this last problem has solution only for a finite number of pairs $(a, b)$, and this will be done in two different ways, essentially according to whether $\cL$ is a rational multiple of $\omega_B^{-1}$ or not.

We recall here the statement of theorem~\ref{thm: finiteness result}.
\begin{RecallThm}
Let $B$ be a smooth projective surface, and $\cL$ an ample line bundle on $B$. Consider the projective bundle $\PP(\cL^a \oplus \cL^b \oplus \cO_B)$, with $a \geq b \geq 0$. Then only for a finite number of pairs $(a, b)$ the generic anticanonical hypersurface in $\PP(\cL^a \oplus \cL^b \oplus \cO_B)$ is a Calabi--Yau elliptic fibration over $B$.
\end{RecallThm}

\subsubsection{Step 1}\label{section: step1}

In this first step we use the informations provided by table~\ref{table: weights} to determine when some of the cohomology spaces where the coefficients $\alpha_{ijk}$ of \eqref{equation: cubic equation} lie are a priori zero.

Since $L$ is an ample divisor, there exists a suitable integer $n_0$ such that $nL + K_B$ is ample for any $n \geq n_0$. Fix one such $n_0$ (e.g.~the least one), then $H^0(B, c_1(B) - nL) = 0$ for any $n \geq n_0$ for otherwise $c_1(B) - nL = -(nL + K_B)$ would be effective. In particular, there is an infinite number of pairs $(a, b)$ satisfying $2a - b \geq n_0$ in the octant $a \geq b \geq 0$: the divisor $(2a - b)L + K_B$ is ample, hence by the previous argument
\[H^0(B, (b - 2a)L - K_B) = H^0(B, -((2a - b)L + K_B)) = 0,\]
and so the coefficient of $x^3$ in \eqref{equation: cubic equation} is identically $0$ (cf.~table~\ref{table: weights}). Equation \eqref{equation: cubic equation} looks then like
\[F = \cancel{\alpha_{300} x^3} + \alpha_{210} x^2 y + \alpha_{201} x^2 z + \ldots\]
and so $\pi: X \longrightarrow B$ has a distinguished section, given by
\begin{equation}\label{equation: section}
P \longmapsto (1: 0: 0) \in X_P.
\end{equation}

Observe that there is only a finite number of pairs $(a, b)$ in the octant $a \geq b \geq 0$ such that $2a - b < n_0$. For such pairs the generic anticanonical hypersurface in $Z$ is a genus $1$ fibration, but since the equation $F$ defining the variety is general, it is difficult to see if there are sections or not. However, we can ignore them from now on since they are only a finite number.

\begin{center}
\begin{figure}[h]
\includegraphics[width = 0.5\textwidth]{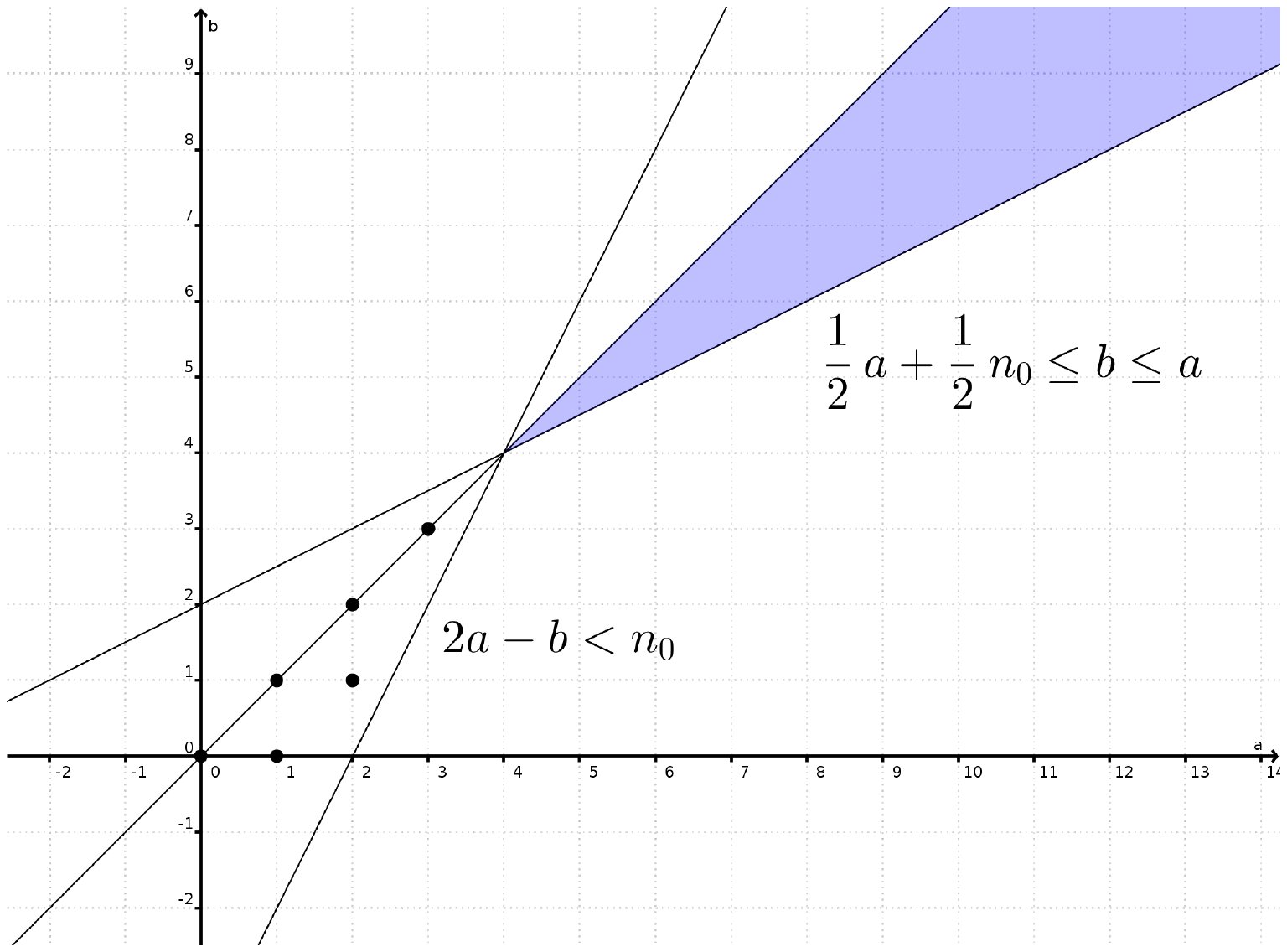}
\caption{The finitely many cases with $2a - b < n_0$. The picture refers to the particular case where $B = \PP^2$ and $L$ is the class of a line, so that $n_0 = 4$. The shaded area corresponds to the bounds given in \eqref{eq: bounds}.}
\label{figure: step1}
\end{figure}
\end{center}

In figure~\ref{figure: step1} it is shown this fact in the particular case where $B = \PP^2$ and $\cL = \cO_{\PP^2}(1)$.

\begin{Rmk}\label{rmk: ineffective bounds}
Exploiting this argument and comparing with the first four rows in table~\ref{table: weights}, it is then easy to see that if $(a, b)$ satisfy
\begin{equation}\label{eq: bounds}
\left\{ \begin{array}{l}
2a - b \geq n_0\\
a \geq n_0\\
b \geq n_0\\
2b - a \geq n_0\\
a \geq b \geq 0
\end{array} \right. \longrightarrow \frac{1}{2} a + \frac{1}{2} n_0 \leq b \leq a,
\end{equation}
then the coefficients $\alpha_{ij0}$ are all necessarily identically zero. In particular, equation \eqref{equation: cubic equation} factors as $F(x, y, z) = z \cdot f(x, y, z)$ and so $F = 0$ can not define a smooth variety. However, this is not enough to conclude the proof of our main theorem, since we are still let with an infinite number of pairs $(a, b)$'s.
\end{Rmk}

\subsubsection{Step 2}\label{section: step2}
It follows from the first step that in the infinitely many cases where $2a - b \geq n_0$, the generic anticanonical hypersurface of $Z$ admits the presence of a section, as defined in \eqref{equation: section}. In this step, we want to to use the relation in proposition~\ref{prop: friedman} to drop down the problem to $B$.

Let $S$ be the image of the section \eqref{equation: section}, by proposition~\ref{prop: friedman} we have that
\[c_2(X) [S] = c_2(S) - c_1(S)^2 = c_2(B) - c_1(B)^2,\]
and so we need to compute the term on the left.

Let $i: X \hookrightarrow Z$ be the inclusion: by \eqref{formula: chern classes of X}, we have that $c_2(X) = i^* \psi$, where
\[\psi = 3 \xi^2 + p^*(2(a + b)L + 3 c_1(B)) \xi + p^*((a + b)L c_1(B) + ab L^2 + c_2(B))\]
and so
\[\deg c_2(X)[S] = \deg i^* \psi \cdot [S] = \deg i_*(i^* \psi \cdot [S]) = \deg \psi \cdot i_*[S].\]
In order to compute $i_*[S]$, which is the class of $S$ in the Chow ring of $Z$, we recall that $S$ is defined in $Z$ by $y = z = 0$, and that this intersection is transverse. So
\[i_*[S] = (\xi + b p^* L) \xi = \xi^2 + b p^* L \xi,\]
and the relation $\deg \psi \cdot i_*[S] = c_2(B) - c_1(B)^2$ reduces to
\begin{equation}\label{equation: conic}
a(a - b)L^2 + (b - 2a)c_1(B)L + c_1(B)^2 = 0.
\end{equation}

Observe that now we have a problem concerning only the base and its intersection theoretic properties. Thinking to $(a, b) \in \R^2$, equation \eqref{equation: conic} defines a plane conic, which is reducible if and only if
\[L^2 = 0 \qquad \text{or} \qquad (c_1(B)L)^2 = L^2 c_1(B)^2.\]
The first case is impossible since we are assuming that $L$ is ample.

By the Hodge index theorem, $(c_1(B)L)^2 \geq L^2 c_1(B)^2$ and
\begin{equation}\label{equation: HIT consequence}
(c_1(B)L)^2 = L^2 c_1(B)^2 \Longleftrightarrow r L \equiv s c_1(B)
\end{equation}
for suitable integers $r$ and $s$ (where $\equiv$ denotes numerical equivalence), and $s \neq 0$ since $L$ is ample.

Our next step is to study the conic defined in \eqref{equation: conic} when it is irreducible (section~\ref{section: step 3irr}) and when it is reducible (section~\ref{section: reducible conic case}), and to show that in each of these two cases we have only a finite number of integral points $(a, b)$ in the octant $a \geq b \geq 0$ on the conic \eqref{equation: conic}.

\subsubsection{Step 3, case 1}\label{section: step 3irr}
Let's concentrate first on the case when the conic \eqref{equation: conic} is irreducible: it is a hyperbola, with asymptotes
\[a = \frac{c_1(B)L}{L^2} \qquad \text{and} \qquad b = a - \frac{c_1(B)L}{L^2}.\]

Observe that if we multiply \eqref{equation: conic} by $L^2$, then it can be written as
\[(L^2 a - c_1(B)L)(L^2 (a - b) - c_1(B)L) = (c_1(B)L)^2 - c_1(B)^2 L^2\]
and so the integral points of \eqref{equation: conic} are the integral pairs $(a_i, b_i)$ having
\[a_i = \frac{d_i + c_1(B)L}{L^2}, \qquad b_i = \frac{d_i - d'_i}{L^2} = \frac{d_1^2 + c_1(B)^2 L^2 - (c_1(B) L)^2}{L^2 d_i},\]
where $d_i$ runs through all the divisors of $(c_1(B)L)^2 - c_1(B)^2 L^2$, and $d'_i = \frac{(c_1(B)L)^2 - c_1(B)^2 L^2}{d_i}$.

So it is clear that they are finite.

\subsubsection{Step 3, case 2}\label{section: reducible conic case}
We concentrate now on the case where the conic \eqref{equation: conic} is reducible, i.e.\ the case where $(c_1(B)L)^2 = L^2 c_1(B)^2$.

The equation for the conic \eqref{equation: conic} is
\[(L^2 a - c_1(B)L)(L^2 a - L^2 b - c_1(B)L) = 0.\]
By \eqref{equation: HIT consequence}, $r L \equiv s c_1(B)$ implies $\frac{c_1(B)L}{L^2} = \frac{r}{s}$: we have two further sub-cases according to whether $\frac{r}{s}$ is a positive integer or not.

If $\frac{r}{s} \notin \N$, the two lines
\[a = \frac{c_1(B)L}{L^2} \qquad \text{and} \qquad b = a - \frac{c_1(B)L}{L^2}\]
have no integral points in the octant $a \geq b \geq 0$. This means that we have no new smooth Calabi--Yau fibrations.

If instead $\frac{r}{s} \in \N$, then in the range $a \geq b \geq 0$ we have a finite number of pairs $(a, b)$ on the line $a = \frac{c_1(B)L}{L^2}$, namely $\frac{c_1(B)L}{L^2} + 1 = \frac{r}{s} + 1$, and an infinite number of $(a, b)$ on the line $b = a - \frac{c_1(B)L}{L^2}$. To give a limitation on the number of these last, we look at the coefficients of the first monomials in equation \eqref{equation: cubic equation}, which are listed in table~\ref{table: weights on the line}, and use the integer $n_0$ introduced in section~\ref{section: step1}. It was defined by the property that $nL - c_1(B)$ is ample for any $n \geq n_0$.

\begin{center}
\begin{longtable}{|c|c|}
\caption{Weight of $\alpha_{ij0}$ on the line $b = a - \frac{c_1(B)L}{L^2}$}
\label{table: weights on the line}\\
\hline
\text{Monomial} & \text{Weight of the coefficient}\\
\hline
\hline
$x^3$ & $c_1(B) - \left( b + 2\frac{r}{s} \right)L$\\
\hline
$x^2 y$ & $c_1(B) - \left( b + \frac{r}{s} \right)L$\\
\hline
$x y^2$ & $c_1(B) - b L$\\
\hline
$y^3$ & $c_1(B) - \left( b - \frac{r}{s} \right)L$\\
\hline
\end{longtable}
\end{center}

Arguing as we did in remark~\ref{rmk: ineffective bounds}, now we find a bound: if $b \geq n_0 + \frac{r}{s}$, we have that all the bundles listed in table~\ref{table: weights on the line} are anti-ample. Hence the coefficients of $x^3$, $x^2 y$, $x y^2$ and $y^3$ in \eqref{equation: cubic equation} are necessarily identically zero, and so the equation $F$ for the variety factors as $F(x, y, z) = z \cdot f(x, y, z)$. Then $F = 0$ can't define a smooth variety.

Observe that in this case $z = 0$ defines a divisor whose class is $\xi$, while $f(x, y, z) = 0$ defines a divisor of class $p^*c_1(B) + (a + b) p^* L + 2\xi$, which is neither a Calabi--Yau variety nor an elliptic fibration.

In particular, we have only a finite number of pairs $(a, b)$ on the line $b = a - \frac{c_1(B)L}{L^2}$ such that the generic anticanonical hypersurface in $\PP(\cL^a \oplus \cL^b \oplus \cO_B)$ could define a Calabi--Yau elliptic fibration over $B$, and a limitation is
\[\frac{r}{s} \leq a \leq n_0 + 2\frac{r}{s} - 1, \qquad 0 \leq b \leq n_0 + \frac{r}{s} - 1.\]
We can be even more precise (see also remark~\ref{rmk: natural setting}), since up to numerical equivalence we have $\frac{r}{s} L \equiv c_1(B)$, and so $nL - c_1(B) \equiv \left( n - \frac{r}{s} \right)L$ is ample if $n \geq \frac{r}{s} - 1$: this means that we can choose $n_0 = \frac{r}{s} - 1$, which gives us the limitations
\begin{equation}\label{formula: limitation}
\frac{r}{s} \leq a \leq 3\frac{r}{s}, \qquad 0 \leq b \leq 2\frac{r}{s}.
\end{equation}

\begin{Rmk}
It is interesting to observe that the ``extreme'' case of limitation \eqref{formula: limitation} occur. In fact choosing $(a, b) = \left( 3 \frac{r}{s}, 2 \frac{r}{s} \right)$, from the relation $r L \equiv s c_1(B)$ we get $3 \frac{r}{s} L \equiv 3 c_1(B)$, $2 \frac{r}{s} L \equiv 2 c_1(B)$, and so we are dealing with the projective bundle
\[\PP(\omega_B^{-3} \oplus \omega_{B}^{-2} \oplus \cO_B),\]
where we can find all the Weierstrass models of the elliptic fibrations over $B$ whose total space is a Calabi--Yau manifold (cf.~section~\ref{sect: fund line bund}).
\end{Rmk}

\begin{center}
\begin{figure}[h]
\includegraphics[width = 0.5\textwidth]{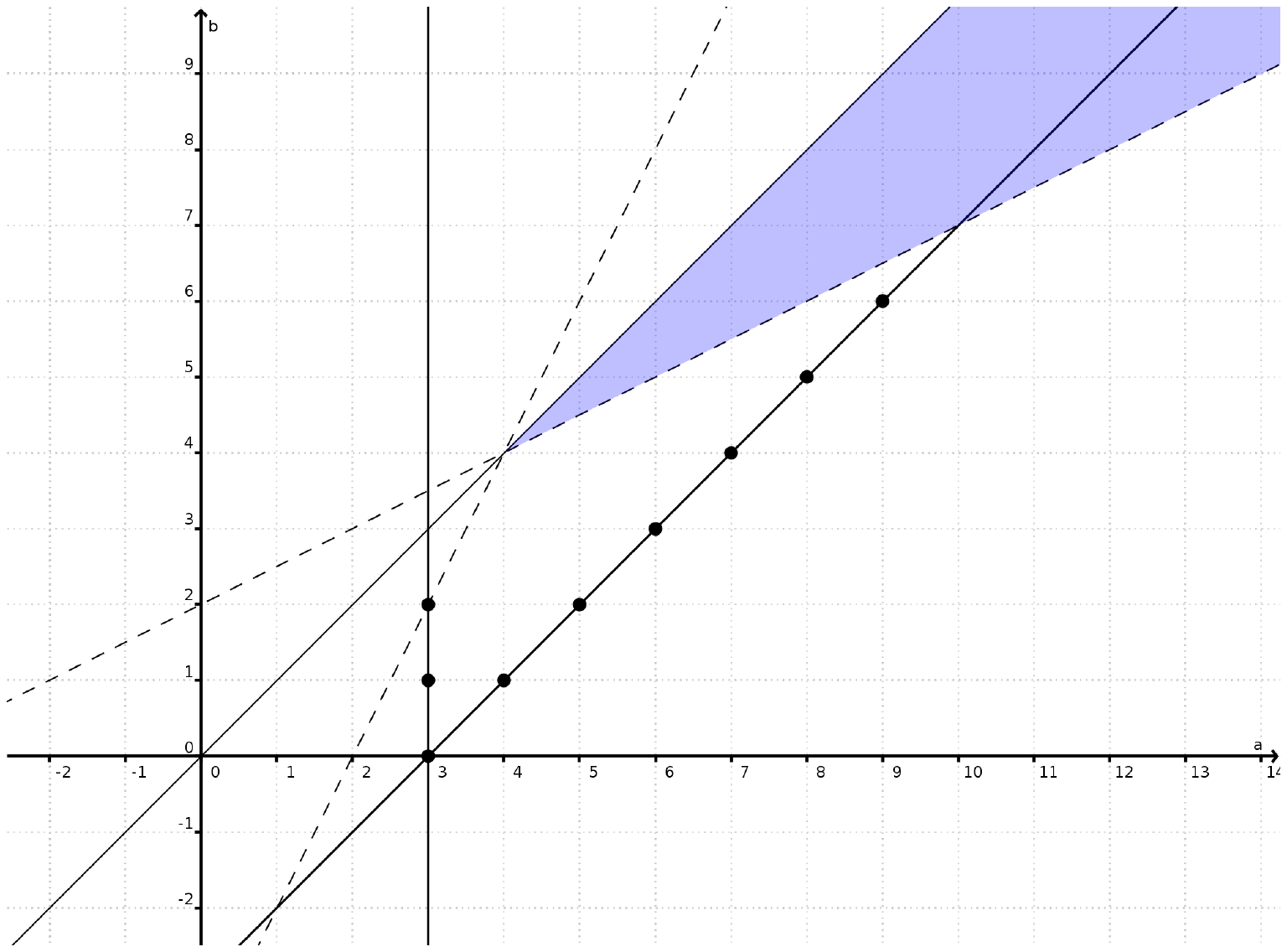}
\caption{If $B = \PP^2$ and $L$ is the class of a line, then we are in the case described in section~ \ref{section: reducible conic case}, and this is the corresponding picture.}
\label{figure: reducible conic}
\end{figure}
\end{center}

If $\frac{r}{s} \in \N$ we have then at most
\begin{equation}\label{formula: estimate}
3 \frac{r}{s} + 1 = \underbrace{\left( \frac{r}{s} + 1 \right)}_{\begin{array}{c}\text{Pairs on the line}\\ a = \frac{c_1(B)L}{L^2}\end{array}} + \underbrace{\left( 2\frac{r}{s} + 1 \right)}_{\begin{array}{c}\text{Pairs on the line}\\ b = a - \frac{c_1(B)L}{L^2}\end{array}} - \underbrace{1}_{\begin{array}{c}\text{The common case}\\ (a, b) = \left( \frac{c_1(B)L}{L^2}, 0 \right)\end{array}}
\end{equation}
such pairs $(a, b)$.

\subsubsection{Conclusion}
Only for a finite number of pairs $(a, b)$ the generic anticanonical hypersurface in $\PP(\cL^a \oplus \cL^b \oplus \cO_B)$ is a smooth Calabi--Yau elliptic fibration, which completes the proof of theorem~\ref{thm: finiteness result}.

We summarize the results obtained in table~\ref{table: summary}.

\begin{center}
\begin{table}[h!]
\caption{Summary of the results}
\label{table: summary}
\begin{tabular}{|p{0.20\textwidth}||p{0.3\textwidth}|p{0.17\textwidth}|p{0.17\textwidth}|}
\hline
\multirow{3}{0.20\textwidth}{$(2a - b)L + K_B$ is not ample} & \multicolumn{3}{c|}{$(2a - b)L + K_B$ is ample}\\
\cline{2-4}
 & \multirow{2}{0.3\textwidth}{$(K_B L)^2 \neq K_B^2 L^2$} & \multicolumn{2}{c|}{$(K_B L)^2 = K_B^2 L^2$}\\
\cline{3-4}
 & & $\frac{r}{s} \notin \N$ & $\frac{r}{s} \in \N$\\
\hline
\hline
Finite number of cases, which are a priori only genus one fibrations. It is not clear if they have at least one section or not. & The conic \eqref{equation: conic} is irreducible, and we have a finite number of Calabi--Yau elliptic fibrations. & No pairs. & Finite number of Calabi--Yau elliptic fibrations, at most $3 \frac{r}{s} + 1$.\\
\hline
\end{tabular}
\end{table}
\end{center}

\begin{Rmk}
We want to stress that we proved that the number of genus $1$ fibrations whose total space is smooth lie in a finite number of $\PP(\cL^a \oplus \cL^b \oplus \cO_B)$, but we don't know a priori if all of them are elliptic fibrations. In the finite number of cases detected in section~\ref{section: step1} it is not clear in fact if there is at least a section.
\end{Rmk}

\begin{Rmk}
We can read our result also in another way: only for a finite number of $Z = \PP(\cL^{\otimes a} \oplus \cL^{\otimes b} \oplus \cO_B)$ the generic element of the anticanonical system $|-K_Z|$ is a \emph{smooth} hypersurface. Let now focus on the infinite number of cases where this does not hold: in view of Bertini's Theorem we can then claim that for such ambient spaces $Z$, the linear system $|-K_Z|$ is not base point free.
\end{Rmk}

\section{Examples}\label{section: examples}

We want to run this program in two cases of interest: the case where the base $B$ is a del Pezzo surface and $L$ is a rational multiple of an anticanonical divisor, and the case where $B$ is a Hirzebruch surface and $L$ is any ample line bundle.

The reason why del Pezzo surfaces are interesting is provided by the following observation.

\begin{Rmk}\label{rmk: natural setting}
Let $B$ be a surface and $L$ an ample divisor on $B$. Assume that at the end of step 2 (section~\ref{section: step2}), the conic \eqref{equation: conic} is reducible. It follows easily from \eqref{equation: HIT consequence} that then $B$ is a del Pezzo surface and $L$ is (numerically) a rational multiple of $c_1(B)$.
\end{Rmk}

Before dealing with the general case in section~\ref{section: del pezzo}, it is worthwhile to study apart the sub-case $B = \PP^2$ (section~\ref{section: over P2}).

The motivation why we are interested in Hirzebruch surfaces if the following.

\begin{Rmk}
Assume that $\pi: X \longrightarrow B$ is a smooth elliptic Calabi--Yau threefold, with $B$ a smooth minimal surface. It follows from \cite[Cor.~3.3]{Grassi} and \cite[Thm.~3.1]{Grassi} that either $B$ is birationally ruled or $B$ is a K3 or Enriques surface and the $j$-invariant function is constant. In the first case, it follows from the discussion after \cite[Cor.~3.3]{Grassi} that $B$ can be either $\PP^2$ or a geometrically ruled surface with Sakai invariant $e$ bounded by $0 \leq e \leq 12$. Finally, by \cite[Main Theorem]{Oguiso} we deduce that $B$ is rational, hence it is $\PP^2$ or a Hirzebruch surface $\F_e$ (with $e \neq 1$).
\end{Rmk}

We will deal with Hirzebruch surfaces in section~\ref{section: hirzebruch}.

\subsection{The case of \texorpdfstring{$B = \PP^2$}{B = P2}}\label{section: over P2}
Observe that if $B$ is a smooth surface with $\Pic B \simeq \Z$, then we are necessarily in the case described in section~\ref{section: reducible conic case}.

Take $B = \PP^2$, and $L = dl$ for $d \in \N$ and $l$ a line in $\PP^2$ (figure~\ref{figure: step1} and \ref{figure: reducible conic} correspond to the choice $d = 1$). Now we compute the least integer $n_0$ such that $n_0 L + K_{\PP^2}$ is ample:
\[n_0 = \left\{ \begin{array}{ll}
4 & \text{if } d = 1\\
2 & \text{if } d = 2, 3\\
1 & \text{if } d \geq 4,
\end{array} \right.\]
so the cases satisfying $2a - b < n_0$ (section~\ref{section: step1}) are
\[\begin{array}{ll}
(0, 0), (1, 0), (1, 1), (2, 1), (2, 2), (3, 3) & \text{if } d = 1\\
(0, 0), (1, 1) & \text{if } d = 2, 3\\
(0, 0) & \text{if } d \geq 4.
\end{array}\]
Since $c_1(\PP^2) = 3l$, we have
\[rdl = 3sl \Longleftrightarrow rd = 3s \Longleftrightarrow \frac{r}{s} = \frac{3}{d}.\]
We have only two cases where the ratio $\frac{r}{s}$ is an integer, which correspond to
\[d = 1 \qquad \text{and} \qquad d = 3,\]
i.e.~$L = l$ or $L = -K_{\PP^2}$. For all the other cases, the only possible pair is then $(a, b) = (0, 0)$, with the exception of $L = 2 l$, which has also $(a, b) = (1, 1)$.

For $d = 3$, there are five possibilities: besides the two we already know, on the reducible conic \eqref{equation: conic} we have also the pairs $(a, b) = (1, 1), (2, 1), (2, 3)$.

\begin{center}
\begin{longtable}{|c||c|}
\caption{Summary of cases with $B = \PP^2$, $L = dl$ and $d \geq 2$.}\\
\hline
$d$ & Possible $(a, b)$\\
\hline
\hline
$2$ & $(0, 0)$, $(1, 1)$\\
\hline
$3$ & $(0, 0)$, $(1, 0)$, $(1, 1)$, $(2, 1)$, $(2, 3)$\\
\hline
$\geq 4$ & $(0, 0)$\\
\hline
\end{longtable}
\end{center}

The only case left is $d = 1$ in the situation of section~\ref{section: reducible conic case}. We have to count the integral points on the conic
\[(a - 3)(a - b - 3) = 0\]
which are in the first octant and have $b \leq 6$ (estimate \eqref{formula: limitation}). On the line $a = 3$ we have the points $(3, 2)$, $(3, 1)$ and $(3, 0)$, while on the line $b = a - 3$ the points $(4, 1)$, $(5, 2)$, $(6, 3)$, $(7, 4)$, $(8, 5)$ and $(9, 6)$.

Then the pairs $(a, b)$ such that the generic anticanonical hypersurface in the bundle $\PP(\cO_{\PP^2}(a) \oplus \cO_{\PP^2}(b) \oplus \cO_{\PP^2})$ could be a smooth Calabi--Yau elliptic fibration are the following $15$:
\[\begin{array}{c}
(0, 0), \qquad (1, 0), \qquad (1, 1), \qquad (2, 1), \qquad (2, 2), \qquad (3, 3),\\
(3, 2), \qquad (3, 1), \qquad (3, 0),\\
(4, 1), \qquad (5, 2), \qquad (6, 3), \qquad (7, 4), \qquad (8, 5), \qquad (9, 6).
\end{array}\]

\begin{Rmk}
Some of these families are already known. For example, the families corresponding to $(a, b) = (3, 3)$ and $(6, 3)$ were analysed in \cite{AluffiEsole}, while the one corresponding to $(a, b) = (6, 3)$ and $(3, 0)$ were analysed in \cite{CacciatoriCattaneoVanGeemen}.
\end{Rmk}

\subsection{The case of del Pezzo surfaces}\label{section: del pezzo}
Let $B$ denote a del Pezzo surface and $\cL$ a rational multiple of the anticanonical bundle, say $\cL^r = \omega_B^{-s}$ (this is the natural setting by remark~\ref{rmk: natural setting}). Let $n_0 = \left[ \frac{r}{s} \right] + 1$, then $nL + K_B$ is ample for all $n \geq n_0$. With the notation of section~\ref{section: step1}, the number of pairs $(a, b)$ for which we can not ensure the presence of a section, i.e.~those satisfying the system
\[\left\{ \begin{array}{l}
a \geq b \geq 0\\
2a - b < n_0,
\end{array} \right.\]
is
\begin{equation}\label{formula: explicit bound 1}
\frac{n_0 (n_0 + 2)}{4} \text{ for $n_0$ even}, \qquad \frac{n_0^2 + 4 n_0 - 1}{4} \text{ for $n_0$ odd}.
\end{equation}

If the ratio $\frac{r}{s}$ is not an integer, then these are the only cases among which we can find elliptic fibrations.

\begin{Rmk}
In particular, for $r < s$ we have only the pair $(a, b) = (0, 0)$.
\end{Rmk}

If the ratio $\frac{r}{s}$ is an integer $m$, then $r = m s$ and so $mL = -K_B$, i.e.~$L$ is a submultiple of $-K_B$. In this case $n_0 = m + 1$ and we have to count also the points on the reducible conic \eqref{equation: conic}: in view of estimate \eqref{formula: estimate} these are $3m$ since the point $(a, b) = (m, m)$ was already taken into account. But then the number of families of elliptic Calabi--Yau threefolds over $B$ is bounded by
\begin{equation}\label{formula: explicit bound 2}
\frac{m^2 + 18 m + 4}{4} \text{ for $m$ even}, \qquad \frac{m^2 + 16 m + 3}{4} \text{ for $m$ odd}.
\end{equation}

\begin{Rmk}
Observe that these results agree with the ones we found in section~\ref{section: over P2} for the plane $\PP^2$. Let $l$ be the class of a line, then:
\begin{enumerate}
\item For $L = l$, we have $r = 3$, $s = 1$ and so we can use \eqref{formula: explicit bound 2} with $m = 3$: we have $15$ cases.
\item For $L = 2l$, we have $r = 3$, $s = 2$ and so we can use \eqref{formula: explicit bound 1} with $n_0 = 2$: we have $2$ cases.
\item For $L = 3l$, we have $r = s = 1$ and so we can use \eqref{formula: explicit bound 2} with $m = 1$: we have $5$ cases.
\item For $L = kl$, with $k \geq 4$, we have $\frac{r}{s} < 1$ and so we can use \eqref{formula: explicit bound 1} with $n_0 = 1$: we have only one case.
\end{enumerate}
\end{Rmk}

\subsection{The case of Hirzebruch surfaces}\label{section: hirzebruch}
Let $\F_e = \PP(\cO_{\PP^1}(e) \oplus \cO_{\PP^1})$ be a Hirzebruch surface. Then the Picard group of $\F_e$ is generated by two classes, $C$ and $f$, with $C^2 = -e$, $C \cdot f = 1$ and $f^2 = 0$. The canonical divisor of $\F_e$ is $K_{\F_e} = -2 C - (e + 2) f$, and a divisor $L = \alpha C + \beta f$ is ample if and only if (cf.~\cite[Cor.~V.2.18]{HAG})
\begin{equation}\label{eq: ample condition}
\left\{ \begin{array}{l}
\alpha > 0\\
\beta > \alpha e.
\end{array} \right.
\end{equation}

It is then easy to see that $-K_{\F_e}$ is ample if and only if $e < 2$, and so the only minimal Hirzebruch surface which is also a del Pezzo surface is $\F_0 = \PP^1 \times \PP^1$. In what follows, we will then assume that $e \geq 2$.

Following the lines of the proof of theorem~\ref{thm: finiteness result}, we first compute the less integer $n_0$ such that $K_{\F_e} + n L$ is ample for every $n \geq n_0$. Thanks to \eqref{eq: ample condition}, we have that
\[n_0 = \left\{ \begin{array}{ll}
3 & \text{if } \alpha = 1\\
2 & \text{if } \alpha = 2\\
1 & \text{if } \alpha \geq 3.
\end{array} \right.\]
So, as a first result, the pairs $(a, b)$ satisfying $2a - b < n_0$ are
\[\begin{array}{ll}
(0, 0), (1, 0), (1, 1), (2, 2) & \text{if } \alpha = 1\\
(0, 0), (1, 1) & \text{if } \alpha = 2\\
(0, 0) & \text{if } \alpha \geq 3.
\end{array}\]

Next, we have to consider the conic \eqref{equation: conic} and find its integral points in the octant $a \geq b \geq 0$. Observe that \eqref{equation: conic} can be written also as
\[(K_{\F_e} + aL)^2 = bL(K_{\F_e} + aL),\]
which is easier to deal with. Before we find the integral points on this conic, we make a little digression, giving some useful estimate for some intersection numbers. We have that
\[-\frac{K_{\F_e} \cdot L}{L^2} = \frac{2\beta - e\alpha - 2\alpha}{\alpha(2\beta - e\alpha)} = \frac{1}{\alpha} + \frac{2}{2\beta - e\alpha}.\]
Observe then that:
\begin{enumerate}
\item If $\alpha = 1$, we have $\beta > e \geq 2$ and so $\beta \geq 3$. As a consequence, $2\beta - e = \beta + (\beta - e) \geq 3 + 1 = 4$, so we deduce that
\[1 < -\frac{K_{\F_e} \cdot L}{L^2} \leq 1 + \frac{1}{2} = \frac{3}{2}.\]
\item If $\alpha \geq 2$, arguing as above we deduce that $\beta \geq 5$ and so $2\beta - e\alpha = \beta + (\beta - e\alpha) \geq 5 + 1 = 6$. This means that
\[0 < -\frac{K_{\F_e} \cdot L}{L^2} \leq \frac{1}{2} + \frac{1}{3} = \frac{5}{6}.\]
\end{enumerate}

With this estimates, we can then prove the following lemma.

\begin{Lemma}\label{lemma: no int pts}
Let $L$ be an ample line bundle on the Hirzebruch surface $\F_e$ with $e \geq 2$. Then the conic \eqref{equation: conic} has no integral points $(a, b)$ with $a \geq 3$.
\end{Lemma}
\begin{proof}
Write $L = \alpha C + \beta f$ as above, and observe that that we have $-K_{\F_e} \cdot L \geq 6$: this means in particular that \emph{the intersection of any ample divisor with the canonical divisor is strictly negative}. We split the proof in two parts, according to whether $\alpha = 1$ or $\alpha \geq 2$.

If $\alpha = 1$, the oblique asymptote $b = a + \frac{K_{\F_e}}{L^2}$ has $-\frac{3}{2} \leq \frac{K_{\F_e}}{L^2} < -1$ and so it suffices to show that given any integer $a \geq 3$, the $b$-coordinate of the point $(a, b)$ on the conic \eqref{equation: conic} satisfies the inequality $b > a - 2$: this means that this point is in between the asymptote and the closest integral point below it, so this point can not be integral. Since we can write our conic as
\[b = \frac{(K_{\F_e} + aL)^2}{L(K_{\F_e} + aL)}\]
and for $a \geq 3$ we have that $(K_{\F_e} + aL)$ is ample, we see that $b > a - 2$ is equivalent to
\[(K_{\F_e} + aL)(K_{\F_e} + 2L) > 0.\]
But writing down this product explicitly, we find that it turns out to be $(2\beta - e - 2)(a - 2)$, which is positive since $L$ is ample and $a \geq 3$. So we are done in this case.

To deal with the case $\alpha \geq 2$, we argue in the same way, but since we have $-\frac{5}{6} \leq \frac{K_{\F_e} \cdot L}{L^2} < 0$ we want to show that
\[b = \frac{(K_{\F_2} + aL)^2}{L(K_{\F_e} + aL)} > a - 1.\]
This is equivalent to $(K_{\F_e} + aL)(K_{\F_e} + L) > 0$, which is true if $\alpha \geq 3$ since it is the intersection of two ample divisors. It remains to show that the inequality holds also when $\alpha = 2$. Writing explicitly the intersection product, we find that $(K_{\F_e} + aL)(K_{\F_e} + L) = 2(a - 1)(\beta - e - 2)$, which is positive since $L$ is ample and $a \geq 3$.
\end{proof}

Thanks to the previous lemma, it remains to deal with only $5$ integral points in the plane.
\begin{enumerate}
\item The point $(2, 0)$ belongs to the conic \eqref{equation: conic} if and only if $(K_{\F_e} + 2L)^2 = 0$. So, if $\alpha \geq 2$ it can not be a point of the conic, as $K_{\F_e} + 2L$ is ample. On the contrary, if $\alpha = 1$, then $K_{\F_e} + 2L = (2\beta - e - 2)f$ and so $(K_{\F_e} + 2L)^2 = 0$, so we do have an integral point.
\item The point $(2, 1)$ belongs to the conic \eqref{equation: conic} if and only if $(K_{\F_e} + 2L)(K_{\F_2} + L) = 0$. We can assume that $\alpha \neq 1$ since we know that in this case the conic passes through the point $(2, 0)$. We can also discard all the cases with $\alpha \geq 3$, since the intersection on the right is the intersection of two ample divisors. So we are left only with the case where $\alpha = 2$, in which case we have $(K_{\F_e} + 2L)(K_{\F_2} + L) = 2(\beta - e - 2) = 0$. The only possible line bundle is then $L = 2 C + (e - 2) f = -K_{\F_e}$, but we must discard this possibility as $-K_{\F_e}$ is not ample.
\item The point $(2, 2)$ belongs to the conic \eqref{equation: conic} if and only if $K_{\F_e} \cdot (K_{\F_e} + 2L) = 0$. As before, we can assume that $\alpha \geq 2$, in which case $K_{\F_e} + 2L$ is ample. But then, as pointed out in the proof of lemma~\ref{lemma: no int pts}, its intersection with an anticanonical divisor is negative, hence we do not have new integral points.
\item The point $(1, 0)$ belongs to the conic \eqref{equation: conic} if and only if $(K_{\F_e} + L)^2 = 0$. Recall that the points we are considering now must also satisfy $b \leq 2a - n_0$, so we can assume that $\alpha \geq 2$. As before, if $\alpha \geq 3$ we have the self-intersection of an ample divisor, so it can not be zero. If $\alpha = 2$, then $K_{\F_e} + L = (\beta - e - 2) f$ and so $(K_{\F_e} + L)^2 = 0$. This means that we have an integral point.
\item The point $(1, 1)$ belongs to the conic \eqref{equation: conic} if and only if $K_{\F_e} \cdot (K_{\F_e} + L) = 0$. Because of the limitation $b \leq 2a - n_0$, we can restrict to $\alpha \geq 3$. In this case $K_{\F_e} + L$ is ample and so its intersection with an anticanonical divisor is negative, hence we do not have new integral points.
\end{enumerate}

We can then sum up these results in the following proposition.

\begin{Prop}
Let $\F_e$ be a Hirzebruch surface, with $e \geq 2$. Let $L = \alpha C + \beta f$ be an ample divisor on $\F_e$, corresponding to the line bundle $\cL$. Then the generic anticanonical divisor in $\PP(\cL^a \oplus \cL^b \oplus \cO_{\F_e})$ defines a smooth Calabi--Yau elliptic fibration over $\F_e$ only if $(a, b)$ is one in the following list:
\[\begin{array}{ll}
(0, 0), (1, 0), (1, 1), (2, 2); (2, 0) & \text{if } \alpha = 1,\\
(0, 0), (1, 1); (1, 0) & \text{if } \alpha = 2,\\
(0, 0) & \text{if } \alpha \geq 3.
\end{array}\]
\end{Prop}

\begin{Rmk}
Concerning the surface $\F_0 = \PP^1 \times \PP^1$, we have that $\Pic \F_0$ is generated by two classes, $f_1$ and $f_2$, with intersections $f_1^2 = f_2^2 = 0$, $f_1 \cdot f_2 = 1$. The canonical divisor is $K_{\F_0} = -2 f_1 - 2 f_2$, and $-K_{\F_0}$ is ample. So, we have to distinguish two cases.
\begin{enumerate}
\item The line bundle $L$ is $L = f_1 + f_2$. In this case, we can apply the arguments of section~\ref{section: del pezzo}, and we see that the possible pairs $(a, b)$ are
\[\begin{array}{c}
(0, 0), (1, 0), (1, 1), (2, 2),\\
(2, 0), (2, 1),\\
(3, 1), (4, 2), (5, 3), (6, 4).
\end{array}\]
\item The line bundle $L = \alpha f_1 + \beta f_2$ is not a rational multiple of $-K_{\F_0}$. In this case, up to switch $f_1$ and $f_2$, it is not restrictive to assume that $\beta > \alpha > 0$, and arguing as we did in this section we can conclude that the possible pairs $(a, b)$ are the following:
\[\begin{array}{ll}
(0, 0), (1, 0), (1, 1), (2, 2); (2, 0) & \text{if } \alpha = 1,\\
(0, 0), (1, 1); (1, 0) & \text{if } \alpha = 2,\\
(0, 0) & \text{if } \alpha \geq 3.
\end{array}\]
\end{enumerate}
\end{Rmk}

\bibliographystyle{plain}
\bibliography{BiblioOverSurfaces}

\end{document}